\documentclass{amsart}
\usepackage[latin1]{inputenc}
\newtheorem{thm}{Theorem}

\newtheorem{lem}[thm]{Lemma}
\newtheorem{prop}[thm]{Proposition}

\theoremstyle{definition}
\newtheorem{defn}[thm]{Definition}

 \numberwithin{equation}{section}
\newcommand{\To}{\longrightarrow}

\begin{document}

\thanks{Author supported by FPU grant and grant MTM2005-08379 from
MEC of Spain and grant 00690/PI/2004 of Fundanción Séneca of
Región de Murcia.} 
\title[]{Renormings of the dual of James tree spaces}%
\author{Antonio Avilés}
\address{Departamento de Matemáticas\\ Universidad de Murcia\\ 30100 Espinardo (Murcia)\\ Spain}

\begin{abstract}
We discuss renorming properties of the dual of a James tree space $JT$. We present examples of weakly Lindel\"of
determined $JT$ such that $JT^\ast$ admits neither strictly convex nor Kadec renorming and of weakly compactly
generated $JT$ such that $JT^\ast$ does
not admit Kadec renorming although it is strictly convexifiable.
\end{abstract}

\maketitle

The norm of a Banach space is said to be locally uniformly rotund
(LUR) if for every for every $x_0$ with $\|x_0\|=1$ and every
$\varepsilon>0$ there exists $\delta>0$ such that
$\|x-x_0\|<\varepsilon$ whenever $\|\frac{x+x_0}{2}\|>1-\delta$. A
lot of research during the last decades has been devoted to
understanding which Banach spaces have an equivalent LUR norm, and
this is still a rather active line of research. In this note we
are concerned with this problem in the case of dual Banach spaces.
It is a consequence of a result of Fabian and Godefroy
\cite{FabGodPRI} that the dual of every Asplund Banach space (that
is, a Banach space such that every separable subspace has a
separable dual) admits an equivalent norm which is locally
uniformly rotund. It is natural to ask whether, more generally,
the dual of every Banach space not containing $\ell_1$ admits an
equivalent LUR norm. We shall give counterexamples to this
question by looking at the dual of James tree spaces $JT$ over
different trees $T$. However all these examples are nonseparable,
and the problem remains open for the separable case. It was
established by Troyanski~\cite{TroyanskiLUR} that a Banach space
admits an equivalent LUR norm if and only if it admits an
equivalent strictly convex norm and also an equivalent Kadec norm.
We recall that a norm is strictly convex if its sphere does not
contain any proper segment and it is a Kadec norm if the weak and
the norm topologies
coincide on its sphere.\\

In Section  \ref{sectionJT} we shall recall the definition of the
spaces $JT$ and the main properties that we shall need.\\

In Section~\ref{sectionFurther} we remark that the space $JT^\ast$
has a LUR renorming whenever $JT$ is separable, so they cannot
provide any counterexample for the separable case. We also point
out the relation which exists between the renorming properties of
$JT^\ast$ and those of $C_0(\bar{T})$, the space of continuous
functions on the completed tree $\bar{T}$ vanishing at $\infty$.
Haydon~\cite{Haydontrees} gave satisfactory characterizations of
those trees $\Upsilon$ for which $C_0(\Upsilon)$ admits LUR,
strictly convex or Kadec equivalent norm. We show that if
$C_0(\bar{T})$ has a LUR (respectively strictly convex) norm then
also does $JT^\ast$, and that, on the contrary, if $JT^\ast$ has
an equivalent Kadec norm,
then so does $C_0(\bar{T})$. We do not know about any of the converses.\\

In Section \ref{sectionWCGJT} we study the case when $JT$ is weakly compactly generated.  The dual of every
weakly compactly generated space is strictly convexifiable, however we shall show that for some trees, $JT$ is
weakly compactly generated but $JT^\ast$ does not admit any equivalent Kadec
norm.\\

In Section \ref{sectionBaire} we provide a sufficient condition on
a tree $T$ in order that $JT^\ast$ does not admit neither a
strictly convex nor a Kadec renorming, namely that it is an
infinitely branchig Baire tree. This is inspired by a construction
of Haydon which can be found in \cite{ArgMerWLD} of the dual of a
weakly Lindel\"of determined Banach space with no equivalent
strictly convex norm (this space contains nevertheless $\ell_1$).
Similar ideas appear also in other Haydon's papers like
\cite{HaydonBaire} and \cite{Haydontrees}. If we consider a
particular tree constructed by Todor\v{c}evi\'c
\cite{Todorcevicorder}, then the Banach space that we construct is
in addition weakly Lindel\"of determined. The short proof of the
properties of the mentioned tree of Todor\v{c}evi\'c presented in
\cite{Todorcevicorder} is based on metamathematical arguments,
while there exists another proof of Haydon~\cite{HaydonBaire}
using
games. We include another proof in Section \ref{treesection}, purely combinatorial.\\

As we mentioned, it is an open question whether the dual of every
separable Banach space $X$ not containing $\ell_1$ admits an
equivalent LUR norm. For such a space $X$, the bidual ball
$B_{X^{\ast\ast}}$ is a separable Rosenthal compact in the
weak$^\ast$ topology (that is, it is a pointwise compact set of
Baire one functions on a Polish space). Hence, the problem is a
particular instance of the more general whether $C(K)$ is LUR
renormable whenever $K$ is a separable Rosenthal compact.
Todor\v{c}evi\'{c} \cite{TodorcevicnoLUR} has recently constructed
a nonseparable Rosenthal compact $K$ such that $C(K)$ is not LUR
renormable, while Haydon, Molt\'{o} and Orihuela \cite{HayMolOri}
have shown that if $K$ is a separable pointwise compact set of
Baire one functions with countably many discontinuities on a
Polish
space, then $C(K)$ is LUR renormable.\\

This research was done while visiting the National Technical University of Athens. We want to express our
gratitude to the Department of Mathematics for its hospitality. Our special thanks go to Spiros Argyros, the
discussion with whom is the origin of the present work.

\section{General properties of James tree spaces}\label{sectionJT}

In this section we shall give the definition and state some well known of James tree spaces. We recall that a
tree is a partially ordered set $(T,\prec)$ such that for every $t\in T$, the set $\{s\in T : s\prec t\}$ is
well ordered by $\prec$. A chain is a subset of $T$ which is totally ordered by $\prec$ and a segment is a chain
$\sigma$ with the extra property that whenever $s\prec t\prec u$ and $s,u\in\sigma$ then $t\in\sigma$. For a
tree $T$ we consider the James tree space $JT$ which is the completion of $c_{00}(T) = \{f\in\mathbb{R}^T :
|supp(f)|<\omega\}$ endowed with the norm

$$\|f\| = \sup \left\{ \left(\sum_{i=1}^n \left(\sum_{t\in \sigma_i}f(t)\right)^2\right)^{\frac{1}{2}}\right\}$$ where the
supremum runs over all finite families of disjoint segments
$\sigma_1,\ldots,\sigma_n$ of the tree $T$. The space $JT$ is
$\ell_2$-saturated, that is, every subspace contains a copy of
$\ell_2$ and in particular $JT$ does not contain $\ell_1$, cf.
\cite{HagOde} and \cite{ArgMersWCG} and
also \cite{Jamestree}.\\

An element $h^\ast\in\mathbb{R}^T$ induces a linear map
$c_{00}(T)\To \mathbb{R}$ given by $h^\ast(x) = \sum_{t\in T}
h^\ast(t)x(t)$. When such a linear map is bounded for the norm of
$JT$, then $h^\ast$ defines an element of the dual space
$JT^\ast$. This is the case when $h^\ast$ is the characteristic
function of a segment $\sigma$ of the tree, $\chi^\ast_\sigma$,
for which we have indeed $\|\chi^\ast_\sigma\| = 1$. Namely, if we
take an element $x\in c_{00}(T)$ of norm less than or equal to one
we will have, taking only the segment $\sigma$ in the definition
of the norm of $JT$, that $|\sum_{i\in \sigma}x(i)|\leq 1$, and
this is the action
of $\chi^\ast_\sigma$ on $x$.\\

\begin{prop}\label{ell2sums}
If $t_1,\ldots,t_n$ are incomparable nodes of the tree $T$ and we have $f_1,\ldots,f_n\in c_{00}(T)$ such that
all the elements on the support of $f_i$ are greater than or equal to $t_i$, shortly
$|f_i|\leq\chi_{[t_i,\infty)}$, then $$\|f_1+\cdots+f_n\| =
\left(\|f_1\|^2+\cdots+\|f_n\|^2\right)^\frac{1}{2}$$
\end{prop}

Proof: Every segment of the tree $T$ intersects at most one of the
segments $[t_i,\infty)$, so the set whose supremum computes the norm of
$f_1+\cdots+f_n$ consists exactly of the numbers of the form
$\left(\sum_1^n\lambda_i^2\right)^\frac{1}{2}$ where each $\lambda_i$ is
one of the numbers whose supremum
computes the norm of $f_i$.$\qed$\\

An antichain is a subset $S$ of $T$ such that every two different elements of $S$ are incomparable.

\begin{defn}
Let $S$ be an antichain of the tree $T$. We define $X_S$ as the subspace of $JT$ generated by all $x\in
c_{00}(T)$ whose support is contained in $[s,\infty)$ for some $s\in S$. For an element $t\in T$, we denote $X_t
= X_{\{t\}}$.
\end{defn}

The properties of the subspaces $X_S$ are the following:\\

\begin{enumerate}
\item $X_S = \left(\bigoplus_{s\in S} X_s\right)_{\ell_2}$. This
is Proposition \ref{ell2sums}.\\

\item $X_S$ is a complemented subspace of $JT$, indeed we have a
norm one projection $\pi_S: JT\To X_S$ which is defined for an
element $x\in c_{00}(T)$ setting $\pi_S(x)_t = x_t$ if $t\succeq
s$ for some $s\in S$ and $\pi_S(x)_t=0$ otherwise. First, $\pi_S$
reduces the norm because if we have a family of segments providing
a sum for computing the norm of $\pi_S(x)$, then we can assume
that every segment is contained in some $[s,\infty)$ for $s\in S$,
and then, the same segments will provide the same sum for the
computation of the norm of $x$.
Second, clearly $\pi_S(x)=x$ if $x\in X_S$.\\

\item The dual map of the operator $\pi_S$ defined above allows us
to consider $X_S^\ast$ as a subspace of $JT^\ast$ since
$\pi_S^\ast:X_S^\ast\To JT^\ast$ is an isometric embedding because
$\pi_S$ is a projection of norm one. In this way $X_S^\ast$ is
identified with the range of $\pi_S^\ast$, which equals the set of
all elements of $JT^\ast$ which take the same values on $x$ and on
$\pi_S(x)$ for every $x\in JT$ (in particular,
$\chi_{[t,u)}^\ast\in X_s^\ast$ whenever $s\preceq t$). Again,
$X_S^\ast$ is a complemented subspace, since if we call
$i_S:X_S\To JT$ to the inclusion, then $i_S^\ast:JT^\ast\To
X_S^\ast$ is a projection of norm one. Taking duals in (1), we
obtain

$$X_S^\ast= \left(\bigoplus_{s\in S} X^\ast_s\right)_{\ell_2}$$

\item Taking duals again, we have an isometric embedding
$i_S^{\ast\ast}:X_S^{\ast\ast}\To JT^{\ast\ast}$ and a projection
of norm one $\pi_S^{\ast\ast}:JT^{\ast\ast}\To X_S^{\ast\ast}$ and
again
$$X_S^{\ast\ast}= \left(\bigoplus_{s\in S}
X^{\ast\ast}_s\right)_{\ell_2}$$

\end{enumerate}

\section{The relation with $C_0(\bar{T})$}\label{sectionFurther}

We notice first that James tree spaces $JT$ cannot be used to
provide examples of separable Banach spaces with non LUR
renormable dual. Let us denote by $\bar{T}$, the completed tree of
$T$, the tree whose nodes are the initial segments of the tree $T$
(that is, the segments $\sigma$ of $T$ with the property that
whenever $s\prec t$ and $t\in\sigma$ then $s\in\sigma$) ordered by
inclusion. We view $T\subset\bar{T}$ by identifying every $t\in T$
with the initial segment $\{s\in T: s\preceq t\}$. A result of
Brackebusch states that for every tree $T$, $JT^{\ast\ast}$ is
isometric to $J\bar{T}$ where $\bar{T}$ is the completed tree of
$T$. We shall need also that by \cite[Theorem VII.2.7]{DGZ}, if
$Y^\ast$ is a subspace of a weakly compactly generated space,
then $Y$ has an equivalent LUR norm.\\

\begin{prop}
Let $T$ be a tree and $X$ be a separable subspace of $JT$, then
$X^{\ast\ast}$ is a subspace of a weakly compactly generated and
hence, $X^\ast$
admits an equivalent LUR norm.\\
\end{prop}

PROOF: Let $T_1$ be a countable set (that we view as a subtree of
the tree $T$) such that $X\subset\overline{span}(\{\chi_{\{t\}} :
t\in T_1\})\cong JT_1$. Since $T_1$ is a countable tree, it has
countable height $ht(T_1)=\alpha<\omega_1$ and the height of the
completed tree cannot be essentially larger,
$ht(\bar{T}_1)\leq\alpha+1<\omega_1$, so in particular,
$\bar{T}_1$ is countable union of antichains and $J\bar{T}_1$ is
weakly compactly generated. Finally, $JT_1^{\ast\ast} \cong
J\bar{T}_1$, so $X^{\ast\ast}$ is a subspace of a weakly compactly
generated space and so $X^\ast$ is
LUR renormable.$\qed$\\

Let us recall now how Brackebusch identifies the basic elements of
$J\bar{T}$ inside $JT^{\ast\ast}$ in order to get an isometry. For
every initial segment of the tree $T$, $s\in\bar{T}$, we have the
basic element $e_s\in JT^{\ast\ast}$ whose action on every
$x^\ast\in JT^\ast$ is given by:

$$ (\star)\ \ e_s(x^\ast) = \lim_{t\in s}x^\ast(\chi_{\{t\}}).$$

The initial segment $s$ is well ordered, so when we write $a =
\lim_{t\in s}a_t$ we mean that for every neighborhood $U$ of $a$
there exists $t_0\in s$ such that $a_t\in U$ whenever $t\geq t_0$.
We consider a tree $\Upsilon$ endowed with its natural locally
compact topology with intervals of the form $(s,s']$ as basic open
sets, and $\Upsilon\cup\{\infty\}$ its one-point compactification.
Let us
first notice the following fact:\\

\begin{prop}\label{topology of the tree} The set $\{e_s : s\in\bar{T}\}\cup\{0\}$
is homeomorphic in the weak$^\ast$ topology of $JT^{\ast\ast}$ to
the space $\bar{T}\cup\{\infty\}$ through the natural
correspondence.\\
\end{prop}

Proof: Since $\bar{T}\cup\{\infty\}$ is compact, it is enough to
check that the natural identification $\bar{T}\cup\{\infty\}\To
\{e_s : s\in\bar{T}\}\cup\{0\}$ is continuous. The fact that it is
continuous at the points $t\in\bar{T}$ follows immediately from
$(\star)$. For the continuity at $\infty$ we take $V$ a
neighborhood of $0$ in the weak$^\ast$ topology and we shall see
that the set $L = \{t\in\bar{T} : t\not\in V\}$ is a relatively
compact subset of $\bar{T}$. We shall prove that every transfinite
sequence $\{t_\alpha : \alpha<\lambda\}$ of elements of $L$ has a
cofinal subsequence which converges to a point of $L$ (this is a
stronger principle than that every net has a convergent subnet and
holds on those sets with scattered compact closure). A partition
principle due to Dushnik and Miller \cite[Theorem 5.22]{DusMil}
yields that either there is an infinite subsequence
$\{t_{\alpha_n} : n\in\omega\}$ of incomparable elements or there
is a cofinal subsequence in which every couple of elements is
comparable. The first possibility is excluded because we know that
a family of vectors of $J\bar{T}$ corresponding to an antichain is
isometric to the basis of $\ell_2$, and in particular it weakly
(and hence weak$^\ast$) converges to 0, contradicting that $V$ is
a weak$^\ast$ neighborhood of 0. In the latter case, the cofinal
subsequence is contained in a branch of the tree which is a well
ordered set, and again the same partition principle of Dushkin and
Miller \cite[Theorem 5.22]{DusMil} implies that it has a further
cofinal and increasing subsequence, and this
subsequence converges to its lowest upper bound in $\bar{T}$.$\qed$\\

Proposition~\ref{topology of the tree} allows us to view every
element $x^\ast\in JT^\ast$ as a continuous function on
$\bar{T}\cup\{\infty\}$ vanishing at $\infty$, and thus to define
an operator,

$$F:JT^\ast\To C_0(\bar{T}).$$

Recall that $C_0(\bar{T})$ stands for the space of real valued
continuous functions on $\bar{T}\cup\{\infty\}$ vanishing at
$\infty$, endowed with the supremum norm $\|\cdot\|_\infty$.
Haydon~\cite{Haydontrees} has characterized the classes of trees
$\Upsilon$ for which the space $C_0(\Upsilon)$ admits equivalent
LUR, Kadec or strictly convex norms. Notice that $F$ is an
operator of norm 1, since $\|F(x^\ast)\|_\infty =
\sup\{|e_s(x^\ast)| : s\in\bar{T}\}\leq \|x^\ast\|$.\\

\begin{thm}\label{relationwithHaydon}
Let $T$ be a tree.\begin{enumerate} \item If $C_0(\bar{T})$ admits
an equivalent strictly convex norm, then $JT^\ast$ also admits an
equivalent strictly convex norm. \item If $C_0(\bar{T})$ admits an
equivalent LUR norm, then $JT^\ast$ also admits an equivalent LUR
norm. \item If $JT^\ast$ admits an equivalent Kadec norm, then
$C_0(\bar{T})$ also admits an equivalent Kadec norm.
\end{enumerate}
\end{thm}

PROOF: Part (1) follows from the fact that $F$ is a one-to-one
operator and one-to-one operators transfer strictly convex
renorming. Moreover, $F$ has the additional property that the dual
operator $F^\ast: C_0(\bar{T})^\ast\To JT^{\ast\ast}\cong
J\bar{T}$ has dense range, because for every dirac measure
$\delta_s$, $s\in\bar{T}$ we have that $F^\ast(\delta_s) = e_s$.
One to one operators whose dual has dense range transfer LUR
renorming \cite{Memoir}, so this proves part (2). Concerning part
(3), we observe that if $|||\cdot|||$ is an equivalent Kadec norm
on $JT^\ast$ and $\rho:\bar{T}\To \mathbb{R}$ is defined by

$$\rho(s) = \inf\{|||\chi^\ast_\sigma||| : s\subset\sigma\}$$

then $\rho:\bar{T}\To \mathbb{R}$ is an increasing function with
no bad points in the sense of~\cite{Haydontrees}, just by the same
argument as in \cite[Proposition 3.2]{Haydontrees}. Hence, by
\cite[Theorem 6.1]{Haydontrees}, $C_0(\bar{T})$ admits an
equivalent Kadec norm.$\qed$\\

We do not know whether any of the converses of
Theorem~\ref{relationwithHaydon} holds true. Concerning part (3),
no transfer result for Kadec norms is available. In the other two
cases, it would be natural to try to imitate Haydon's arguments in
\cite{Haydontrees} using the function $\rho(s) =
\inf\{|||\chi^\ast_\sigma||| : s\subset\sigma\}$ on $JT^\ast$. But
these arguments rely on the consideration of certain special
functions $f\in C_0(\Upsilon)$ which are not available anymore in
$JT^\ast$ which is a rather smaller space.\\

\section{When JT is weakly compactly generated}\label{sectionWCGJT}

In this section we analyze the case when $JT$ is weakly compactly
generated. This property is characterized in terms of the tree as it is
shown in the following result which can be found in \cite{ArgMerWUR}:

\begin{thm}\label{JTWCG}
For a tree $T$ the following are equivalent
\begin{enumerate}
\item $JT$ is weakly compactly generated. \item $JT$ is weakly countably determined. \item $T$ is the union of
countably many antichains. \item $T = \bigcup_{n<\omega}S_n$ where for every $n<\omega$, $S_n$ contains no
infinite chain.
\end{enumerate}
\end{thm}

A tree is union is the union of countably many antichains if and
only if it is $\mathbb{Q}$-embeddable, cf. \cite[Theorem
9.1]{Todorcevicorder}. It happens that for a tree $T$ satisfying
the conditions of Theorem \ref{JTWCG}, the renorming properties of
$JT^{\ast}$ depend on whether the completed tree $\bar{T}$ is
still the union of countably many antichains.

\begin{thm}\label{JTWCGKadec}
Let $T$ be a tree which is the union of countably many antichains. The
following are equivalent:
\begin{enumerate}
\item $\bar{T}$ is also the union of countably many antichains.
\item $JT^\ast$ admits an equivalent Kadec norm.
\item $JT^\ast$ admits an equivalent LUR norm.
\end{enumerate}
\end{thm}

The dual of every weakly compactly generated space admits always
an equivalent strictly convex norm since, by the
Amir-Lindenstrauss Theorem there is a one-to-one operator into
$c_0(\Gamma)$. Hence, that $(2)$ and $(3)$ are equivalent is a
consequence of the result of Troyanski mentioned in the
introduction. On the other hand, we also mentioned in
Section~\ref{sectionFurther} the result of
Brackebusch~\cite{Brackebusch} that for any tree $T$,
$JT^{\ast\ast}$ is isometric to $J\bar{T}$. Hence, if (1) is
verified, then $JT^{\ast\ast}$ is weakly compactly generated and
it follows then by \cite[Theorem VII.2.7]{DGZ} that $JT^\ast$
admits an equivalent LUR norm. Our goal is therefore to prove that
(2) implies (1) but before passing to this we give an example of a
tree $T_0$ which is the union of countably many antichains but the
completion $\bar{T}_0$ does not share this property, so that after
Theorem \ref{JTWCGKadec} $JT_0$ is a weakly compactly generated
space not containing $\ell_1$ and such that $JT_0^\ast$ does not
admit any equivalent Kadec norm, namely
$$T_0 =\sigma'\mathbb{Q} = \{t\subset \mathbb{Q} : (t,<)\text{ is well ordered and }\max(t)\text{ exists}\},$$
where $t\prec s$ if $t$ is a proper initial segment of $s$. For every rational number $q\in\mathbb{Q}$, the set
$S_q = \{t\in T_0 : \max(t)=q\}$ is an antichain of $T_0$, and $T_0 = \bigcup_{q\in\mathbb{Q}}S_q$. The
completed tree $\bar{T}_0$ can be identified with the following tree:
$$T_1 = \sigma\mathbb{Q} = \{t\subset \mathbb{Q} : (t,<)\text{ is well ordered}\},$$
the identification sending every $t\in T_1$ to the initial segment $\{t'\in T_0 : t'\prec t\}$ of $T_0$. The
fact that $T_1$ is not countable union of antichains is a well known result due to Kurepa \cite{KurepasigmaQ},
cf. also \cite{Todorcevicorder}. The reason is the following: suppose there existed $f:T_1\To\mathbb{N}$ such
that $f^{-1}(n)$ is an antichain. Then we could construct by recursion a sequence $t_1\prec t_2\prec\cdots$
inside $T_1$ and a sequence of rational numbers $q_1>q_2>\cdots$ such that $q_i> \sup(t_i)$ and $f(t_{n+1}) =
\min\{f(t) : t_n\prec t, \sup(t)<q_n\}$. The consideration of the element $t_\omega = \bigcup_{n<\omega}t_n$
leads to a contradiction.\\

\begin{lem}\label{Kadecenelarbol}
Let $T$ be any tree and suppose that there exists an equivalent Kadec norm
on $JT$, then there exist \begin{itemize} \item [(a)] a countable
partition of $\bar{T}$, $\bar{T}=\bigcup_{n<\omega}T_n$ and \item [(b)] a
function $F:\bar{T}\To 2^T$ which associates to each initial segment
$\sigma\in \bar{T}$ a finite set $F(\sigma)$ of immediate successors of
$\sigma$,
\end{itemize}
such that for every $n<\omega$ and for every infinite chain
$\sigma_1\prec\sigma_2\prec\cdots$ contained in $T_n$ there exists
$k_0<\omega$ such that $F(\sigma_k)\cap \sigma_{k+1}\neq\emptyset$ for
every $k>k_0$.
\end{lem}

Proof: Let $|||\cdot|||$ be an equivalent Kadec norm on $JT^\ast$.\\

\emph{Claim}: For every $\sigma\in \bar{T}$ there exists a natural number
$n_\sigma$ and a finite set $F(\sigma)\subset T$ of immediate successors
of $\sigma$ such that $\left|\
|||\chi_\sigma^\ast|||-|||\chi_{\sigma'}^\ast|||\ \right|\geq
\frac{1}{n_\sigma}$ for every $\sigma'\in \bar{T}$ such that $\sigma\prec
\sigma'$ and
$F(\sigma)\cap\sigma'=\emptyset$.\\

Proof of the claim: Suppose that there existed $\sigma\in\bar{T}$ failing
the claim. Then, we can find recursively a sequence $\{q_n\}$ of different
immediate succesors of $\sigma$ together with a sequence $\{\sigma_n\}$ of
elements of $\bar{T}$ such that $\sigma\cup\{q_n\}\preceq \sigma_n$ and
$$\left|\ |||\chi_\sigma^\ast|||-|||\chi_{\sigma_n}^\ast|||\ \right|<\frac{1}{n}.$$

Now, $\{\sigma'_n=\sigma_n\setminus\sigma\}$ is a sequence of incomparable
segments of $T$, so the sequence $\{\chi_{\sigma'_n}^\ast\}$ is isometric
to the base of $\ell_2$ and in particular it weakly converges to 0. Hence
the sequence $\chi_{\sigma_n}^\ast = \chi_\sigma^\ast +
\chi_{\sigma'_n}^\ast$ weakly converges to $\chi_\sigma^\ast$, however it
does not converge in norm since $\|\chi_{\sigma'_n}^\ast\|=1$ for every
$n$. Finally, since $|||\chi_{\sigma_n}^\ast|||$ converges to
$|||\chi_{\sigma}^\ast|||$ we obtain, after normalizing, a contradiction
with the fact that $|||\cdot|||$ is a Kadec norm.\\

From the claim we get the function $F$ and also the countable
decomposition setting $T_n = \{\sigma\in\bar{T} : n_\sigma = n\}$. Suppose
that we have an increasing sequence $\sigma_1\prec\sigma_2\prec\cdots$
inside $T_n$. We observe that whenever $F(\sigma_k)\cap \sigma_{k+1} =
\emptyset$ we have that $\left|\
|||\chi_{\sigma_k}^\ast|||-|||\chi_{\sigma_{k'}}^\ast|||\ \right|\geq
\frac{1}{n}$ for all $k'>k$. This can happen only for finitely many $k$'s
because $|||\cdot|||$ is an equivalent norm so it is bounded on the unit
sphere of $JT^\ast$.$\qed$\\

Now we assume that $T$ is union of countably many antichains,
$T=\bigcup_{m<\omega}R_m$, and that it verifies the conclusion of Lemma
\ref{Kadecenelarbol} for a decomposition $\bar{T} = \bigcup_{n<\omega}T_n$
and a function $F$, and we shall show that indeed $\bar{T}$ is the union
of countably many antichains. For every $n<\omega$ and every finite subset
$A$ of natural numbers we consider the set $$S_{n,A} =
\left\{\sigma\in\bar{T} : \sigma\in T_n \text{ and } F(\sigma)\subset
\bigcup_{m\in A}R_m\right\}$$

This gives an expression of $\bar{T}$ as countable union $\bar{T} = \bigcup_{n,A}S_{n,A}$. We shall verify that
this expression verifies condition (4) of Theorem~\ref{JTWCG}. Suppose by contradiction that we had an infinite
chain $\sigma_1\prec\sigma_2\prec\cdots$ inside a fixed $S_{n,A}$. First, since $S_{n,A}\subset T_n$ there
exists $k_0$ such that $F(\sigma_k)\cap\sigma_{k+1}\neq\emptyset$ for every $k>k_0$, say $t_k\in
F(\sigma_k)\cap\sigma_{k+1}\subset\bigcup_{m\in A}R_m$. Then $t_1\prec t_2\prec\cdots$ is an infinite chain of
$T$ contained in $\bigcup_{m\in A}R_m$ which is a finite union of antichains. This contradiction finishes the
proof of Theorem~\ref{JTWCGKadec}.

\section{Spaces with no strictly convex nor Kadec
norms}\label{sectionBaire}

In this section we give a criterion on a tree $T$ in order that $JT^\ast$
admits neither a Kadec norm nor a strictly convex norm. We recall that the
downwards closure of a subset $S$ of a tree $T$ is defined as
$$\hat{S} =
\{t\in T : \exists s\in S : t\preceq s\}.$$

\begin{thm}\label{renormJamestreestar}
Let $T$ be a tree verifying the following properties:

\begin{itemize}
\item[(T1)] Every node of $T$ has infinitely many immediate
succesors.

\item[(T2)] For any countable family of antichains $\{S_n :
n<\omega\}$ there exists $t\in T$ such that
$t\not\in\bigcup_{n<\omega}\hat{S}_n$.
\end{itemize}

Then there is neither a strictly convex nor a Kadec equivalent
norm in $JT^\ast$.\\
\end{thm}

Condition (T2) is called \emph{Baire property} of the tree and
condition (T1) is usually expressed saying that $T$ is an
\emph{infinitely branching tree}. An example of a tree satisfying
properties (T1) and (T2) is the tree whose nodes are the countable
subsets of $\omega_1$ with $s\prec t$ if $s$ is an initial segment
of $t$ (property (T2) is proved by constructing a sequence
$t_1\prec t_2\prec\cdots$ with $t_i\not\in \hat{S}_i$ and taking
$t\succ \bigcup t_i$). A refinement of this construction due to
Todor\v{c}evi\'c \cite{Todorcevicorder} produces a tree with the
additional property that all branches are countable, and this
implies that for this tree $JT$ is weakly Lindel\"of determined
\cite{ArgMerWLD}. This is the example
discussed in Section \ref{treesection}.\\

Along the work of Haydon it is possible to find different results
implying that if a tree $\Upsilon$ is an infinitely branching
Baire tree, then $C_0(\Upsilon)$ (or certain spaces which can be
related to it) has no Kadec or strictly convex norm, cf.
\cite{Haydontrees}, \cite{HaydonBaire}. One may be tempted to use
Theorem~\ref{relationwithHaydon} in conjunction with these results
to get Theorem~\ref{renormJamestreestar}. However, there is a
difficulty since these properties (T1) and (T2) on $T$ are not
easily reflected on the completed tree $\bar{T}$. The tree
$\bar{T}$ is never a Baire tree, since the set $M$ of all maximal
elements verifies that $\hat{M}=\bar{T}$, and even if we try to
remove these maximal elements, the hypothesis that $T$ is
infinitely branching is weaker than the hypothesis that $\bar{T}$
is infinitely branching. We shall do it therefore by hand, using
in any case,
similar arguments as in Haydon's proofs.\\

We assume now that $T$ satisfies (T1) and (T2), we fix an
equivalent norm $(JT^\ast,|||\cdot|||)$ and we shall see that this
norm is neither strictly convex nor a Kadec norm.

\begin{lem}\label{lemaclave}
For any node of the tree $t\in T$ and every $\varepsilon>0$ we can
find another node $s\succ t$ and an element $x_s^{\ast\ast}\in
JT^{\ast\ast}$ with $|||x_s^{\ast\ast}|||^\ast=1$ such that
\begin{enumerate}
\item $\left|\sup\{|||\chi_{[0,u]}^\ast||| : u\succeq
s\}-|||\chi_{[0,s]}^\ast|||\right|<\varepsilon$. \item
$x_s^{\ast\ast}(\chi_{[0,u]}^\ast) \geq
|||\chi_{[0,s]}^\ast|||-\varepsilon$ whenever $s\prec u$.
\end{enumerate}
\end{lem}

PROOF: First we take a node $t'\succ t$ such that
$$\left|\sup\{|||\chi_{[0,u]}^\ast||| : u\succeq
t\}-|||\chi_{[0,t']}^\ast|||\right|<\frac{\varepsilon}{2},$$
 and we find
$x^{\ast\ast}\in JT^{\ast\ast}$ with $|||x_s^{\ast\ast}|||^\ast=1$
such that $x^{\ast\ast}(\chi_{[0,t']})=|||\chi_{[0,t']}^\ast|||$.
We consider the set $S$ of all immediate successors of $t'$ in the
tree $T$ which is an infinite antichain. Then, we can consider the
projection

$$\pi_S^{\ast\ast}:JT^{\ast\ast}\To X_S^{\ast\ast}= \left(\bigoplus_{s\in S} X^{\ast\ast}_s\right)_{\ell_2}$$

Since $S$ is infinite, there must exist $s\in S$ such that
$\|\pi_{s}^{\ast\ast}(x^{\ast\ast})\| < \frac{\varepsilon}{2}$.\\

The elements $s\in T$ and $x^{\ast\ast}_s = x^{\ast\ast}$ are the
desired. Namely, for any $u\succeq s$,
$$x^{\ast\ast}(\chi_{[0,u]}^\ast) =
x^{\ast\ast}(\chi_{[0,t']}^\ast) +
x^{\ast\ast}(\chi_{[s,u]}^\ast),$$ and $\chi_{[s,u]}^\ast\in
X_s^\ast$, so $\pi_s^\ast(\chi_{[s,u]}^\ast)=\chi_{[s,u]}^\ast$
and

\begin{eqnarray*}x^{\ast\ast}(\chi_{[0,u]}^\ast) &=&
x^{\ast\ast}(\chi_{[0,t']}^\ast) +
x^{\ast\ast}(\chi_{[s,u]}^\ast)\\
&=& x^{\ast\ast}(\chi_{[0,t']}^\ast) + x^{\ast\ast}(\pi_s^\ast
(\chi_{[s,u]}^\ast))\\ &=& x^{\ast\ast}(\chi_{[0,t']}^\ast) +
\pi_s^{\ast\ast}(x^{\ast\ast})(\chi_{[s,u]}^\ast)\\ &\geq&
x^{\ast\ast}(\chi_{[0,t']}^\ast)-\frac{\varepsilon}{2}\\ &=&
|||\chi_{[0,t']}^\ast||| - \frac{\varepsilon}{2}\\
&\geq& |||\chi_{[0,s]}^\ast||| - \varepsilon.
\end{eqnarray*}

This guarantees in particular that $|||\chi_{[0,s]}^\ast|||\geq
x^{\ast\ast}(\chi^\ast_{[0,s]})\geq
|||\chi_{[0,t']}^\ast|||-\frac{\varepsilon}{2}$. This together
with the property which follows from
the initial choice of $t'$ gives also property (1) in the lemma and finishes the proof.$\qed$\\

We construct by recursion, using Lemma \ref{lemaclave}, a sequence
of maximal antichains of $T$, $\{S_n : n<\omega\}$ which are
increasing (that is for every $t\in S_{n+1}$,
 there exists $s\in S_{n}$ with $s\prec t$) and such that for every $n<\omega$ and for
every $s\in S_n$ there exists an element $x_s^{\ast\ast}\in
JT^{\ast\ast}$ with $|||x_s^{\ast\ast}|||^\ast=1$ such that
\begin{enumerate}
\item $\left|\sup\{|||\chi_{[0,u]}^\ast||| : u\succ
s\}-|||\chi_{[0,s]}^\ast|||\right|<\frac{1}{n}$. \item
$x_s^{\ast\ast}(\chi_{[0,u]}^\ast) =
x_s^{\ast\ast}(\chi_{[0,s]}^\ast) \geq |||\chi_{[0,s]}^\ast|||-\frac{1}{n}$ whenever $s\prec u$.\\
\end{enumerate}

Now, by property (T2), we can pick $t\in
T\setminus\bigcup_{n<\omega}S_n$. We can find for $t$ a sequence
$s_1\prec s_2\prec\cdots\prec t$ with $s_n\in S_n$.\\

For any $t'\succeq t$ and for every $n<\omega$,

$$|||\chi_{[0,s_n]}^\ast|||-\frac{1}{n} \leq x_{s_n}^{\ast\ast}(\chi_{[0,t']}^\ast) \leq
|||\chi_{[0,t']}^\ast||| \leq \sup_{u\succeq
s_n}|||\chi_{[0,u]}^\ast||| \leq
|||\chi_{[0,s_n]}^\ast|||+\frac{1}{n}.$$

This implies that all the successors of $t$ have the same norm
$|||\cdot|||$ equal to the limit of the norms
$|||\chi_{[0,s_n]}^\ast|||$. If we take $t_1$ and $t_2$ two
immediate succesors of $t$, in addition, for every $n<\omega$

$$|||\frac{\chi_{[0,t_1]}^\ast + \chi_{[0,t_2]}^\ast}{2}|||\geq
x^{\ast\ast}_{s_n}\left(\frac{\chi_{[0,t_1]}^\ast +
\chi_{[0,t_2]}^\ast}{2}\right)\geq
|||\chi_{[0,s_n]}^\ast|||-\frac{1}{n}$$

and passing to the limit

$$ |||\frac{\chi_{[0,t_1]}^\ast + \chi_{[0,t_2]}^\ast}{2}|||\geq
|||\chi_{[0,t_1]}^\ast||| = |||\chi_{[0,t_2]}^\ast|||$$

and this shows that $|||\cdot|||$ is not a strictly convex norm.\\

If now we take a sequence of different immediate succesors of $t$,
$\{t_n : n<\omega\}$, then $\chi_{\{t_n\}}^\ast$ is an element of
norm one of $X_{t_n}^\ast$ and since $$X_{\{t_n : n<\omega\}}^\ast
= \left(\bigoplus_{n<\omega} X^\ast_{t_n}\right)_{\ell_2}$$ the
sequence $(\chi_{\{t_n\}}^\ast : n<\omega)$ is isometric to the
base of $\ell_2$ and in particular it is weakly null. Therefore
$\chi_{[0,t_n]}^\ast$ is a sequence in a sphere which weakly
converges to $\chi_{[0,t]}^\ast$ which is in the same sphere.
However $\|\chi_{[0,t_n]}^\ast - \chi_{[0,t]}^\ast\| =
\|\chi_{[t_n,t_n]}^\ast\|=1$ so this sequence does not converge in
norm. This shows that $|||\cdot |||$ is not a Kadec norm.

\section{About a tree of Todor\v{c}evi\'{c}}\label{treesection}

A subset $A$ of $\omega_1$ is called stationary if the
intersection of $A$ with every closed and unbounded subset of
$\omega_1$ is nonempty. We shall fix a set $A$ such that both $A$
and $\omega_1\setminus A$ are stationary. The existence of such a
set follows from a result of Ulam \cite[Theorem
3.2]{Kunencombinatorics}.

\begin{defn}[Todor\v{c}evi\'c]
We define $T$ to be the tree whose nodes are the closed subsets of $\omega_1$ which are contained in $A$ and
whose order relation is
that $s\prec t$ if $s$ is an initial segment of $t$.\\
\end{defn}

First, $T$ has property (T1) because if $t\in T$ and $\eta\in A$
verifies that $\eta>\max(t)$, then $t\cup \{\gamma\}$ is an
immediate successor of $t$ in $T$. On the other hand, $T$ does not
contain any uncountable chain. If $\{t_i\}{i<\omega_1}$ were an
uncountable chain, then $\bigcup_{i<\omega_1}t_i$ is a closed an
unbounded subset of $\omega_1$, so it should intersect
$\omega_1\setminus A$, which is impossible. The difficult point is
in showing that $T$ verifies property (T2).\\

\begin{thm}[Todor\v{c}evi\'c]\label{antichainstree}
For any countable family of antichains $\{S_n : n<\omega\}$ there
exists $t\in T$ such that $t\not\in\bigcup_{n<\omega}\hat{S}_n$.\\
\end{thm}

PROOF: We suppose by contradiction that we have a family of antichains
$\{S_n : n<\omega\}$ which does not verify the statement. We can suppose
without loss of generality that every one of these antichains is a maximal
antichain, and that they are increasing, that is, for every $t\in S_{n+1}$
there exists $s\in S_n$ such that $s\prec t$. What we know is that for
every $t\in T$ we can find
$t'\in\bigcup_{m<\omega}S_m$ such that $t\prec t'$. Moreover, since the antichains are taken maximal and increasing,\\

$(\ast)$ For every natural number $n$ and for every $t\in T$ there
exists $t'\in\bigcup_{m>n}S_m$ such that $t\prec t'$.\\

We construct a family $\{R_\xi : \xi<\omega_1\}$ of subsets of $T$
with the following properties:\\

\begin{enumerate}
\item $R_\xi$ is a countable subset of $\bigcup_{n<\omega} S_n$.
\item $R_\xi\subset R_\zeta$ whenever $\xi<\zeta$.
\item If $\xi$ is a limit ordinal, then $R_\xi =
\bigcup_{\zeta<\xi}R_\zeta$.
\item If we set $\gamma_\xi = \sup\{\max(t) : t\in R_\xi\}$ then
the following are satisfied
\begin{enumerate}
\item $\gamma_\xi<\gamma_\zeta$ whenever $\xi<\zeta$. \item For
every $\xi<\omega_1$, every $t\in R_\xi$, every $n<\omega$ and
every $\eta\in A$ such that $\max(t)<\eta<\gamma_\xi$ there exists
$t'\in R_\xi\cap\cup_{m>n}S_m$ such that $t\cup\{\eta\}\prec t'$.
\item $\gamma_\xi\neq \max(t)$ for every $t\in R_\xi$.\\
\end{enumerate}
\end{enumerate}

These sets are constructed by induction on $\xi$. We set
$R_0=\emptyset$ and we suppose we have constructed $R_\zeta$ for
every $\zeta<\xi$. If $\xi$ is a limit ordinal, then we define
$R_\xi = \bigcup_{\zeta<\xi}R_\zeta$. Notice that then $\gamma_\xi
= \sup\{\gamma_\zeta : \zeta<\xi\}$ and all properties are
immediately verified for $R_\xi$ provided they are verified for
every $\zeta<\xi$.\\

Now, we suppose that $\xi=\zeta +1$. In order that 4(b) is
verified, we will carry out a saturation argument. We will find
$R_\xi$ as the union of a sequence $R_\xi =
\bigcup_{n<\omega}R_\xi^n$.\\

First, we set $R_\xi^0 = R_\zeta$ and $\gamma_\xi^0 =
\gamma_\zeta$. Because we know that property 4(b) is verified by
$R_\zeta$, we have guaranteed property 4(b) in $R_\xi$ when
$\eta<\gamma_\zeta$.\\

In the next step, we take care that 4(b) is verified for every
$t\in R_\xi^0$ and $\eta=\gamma_\zeta$. That is, for every $t\in
R_\xi^0$ and every $n<\omega$ we find, using property $(\ast)$,
$t'_n\in\bigcup_{m>n}S_m$ such that $t\cup\{\gamma_\xi^0\}\prec
t'_n$ and we set $R_\xi^1 = R_\xi^0\cup\{t'_n : t\in R_\xi^0,\
n<\omega\}$ and $\gamma_\xi^1 = \sup\{\max(s)
: s\in R_\xi^1\}$.\\

If we have already defined $R_\xi^{n}$ and $\gamma_\xi^{n} =
\sup\{\max(s) : s\in R_\xi^{n}\}$ then we make sure that property
4(b) will be verified in $R_\xi$ for any $\eta\leq\gamma_\xi^n$,
that is for every
 every $n<\omega$,
every $t\in R_\xi^{n}$ and every $\eta\in(max(t),\gamma_\xi^n]$,
we find, by property $(\ast)$, an element
$t'_{n\eta}\in\bigcup_{m>n}S_m$ such that
$t\cup\{\gamma_\xi^0\}\prec t'_{n\eta}$ and we set $$R_\xi^{n+1} =
R_\xi^n\cup\{t'_{n\eta} : t\in R_\xi^0,\ n<\omega,\
\eta\in(\max(t),\gamma_\xi^n]\}$$ and $\gamma_\xi^{n+1} =
\sup\{\max(s)
: s\in R_\xi^{n+1}\}.$\\

Finally, setting $R_\xi=\bigcup_{n<\omega} R_\xi^n$, we will have
that $\gamma_\xi = \sup_{n<\omega}\gamma_\xi^n$ and the
construction is finished.\\

Now, we will derive a contradiction from the existence of the sets
$R_\xi$. The set $\{\gamma_\xi : \xi<\omega_1\}$ is a closed and
unbounded subset of $\omega_1$, so since $A$ is stationary, there
exists $\xi<\omega_1$ such that $\gamma_\xi\in A$. We will
construct a sequence $t_1\prec t_2\prec\cdots$ of elements of
$R_\xi$ such that $t_n\in \bigcup_{m>n}S_m$ and $\gamma_\xi =
\sup\{\max(t_n) : n<\omega\}$. Such a sequence leads to a
contradiction, because in this case, $t=\bigcup_{n=1}^\infty
t_n\cup\{\gamma_\xi\}$ is a node of the tree with the property
that for every $n$, $t\succ t_n\in S_{m_n}$, $m_n>n$, and this
implies that $t\not\in\bigcup_{n<\omega}\hat{S}_n$. The
construction of the sequence $t_n$ is done inductively as follows.
An increasing sequence of ordinals $\{\eta_i : i<\omega\}$
converging to $\gamma_\xi$ is chosen. If we already defined
$t_{n-1}$, we find $i$ with $max(t_n)<\eta_i$ and we use property
4(b) to find $t_n\in R_\xi\cap\bigcup_{m>n}S_m$ with
$t_{n-1}\cup\{\eta_i\}\prec t_n$.\\

\end{document}